\documentclass[12pt]{article}

\usepackage{times}
\usepackage[centertags]{amsmath}
\usepackage{amssymb}
\usepackage{amsthm}
\usepackage[latin1]{inputenc}
\usepackage[colorlinks=true,urlcolor=red,citecolor=red,linkcolor=red]{hyperref}
\usepackage[english]{babel}
\usepackage{graphicx}
\usepackage{tikz}
\usepackage{qtree}
\usepackage{tikz-qtree}
\usepackage[lined,algonl,boxed]{algorithm2e}
\usetikzlibrary{trees,positioning,shapes} 
\newtheorem{prop}{Proposition}[section]

\newtheorem{Lemma3}[prop]{Lemma}

\theoremstyle{definition} 
\newtheorem{example1}[prop]{Example}
\newtheorem{example2}[prop]{Example}
\newtheorem{example3}[prop]{Example}
\newtheorem{example4}[prop]{Example}
\newtheorem{example44}[prop]{Example}

\newtheorem{remark0}[prop]{Remark}
\newtheorem{remark1}[prop]{Remark}
\newtheorem{remark2}[prop]{Remark}

\newenvironment{acknowledgements}%
    {\cleardoublepage\thispagestyle{empty}\null\vfill\begin{center}%
    \bfseries Acknowledgements\end{center}}%
    {\vfill\null}

\newcounter{lastnote}

\title{Arf good semigroups with fixed conductor}

\author
{Giuseppe Zito\\
\\
\\
}

\date{}

\begin{document}

\baselineskip16pt

\maketitle

\begin{abstract}
In this paper we present a new algorithm to calculate the set of Arf numerical semigroups
with given conductor. Moreover, we extend the previous problem to the class of good semigroups, presenting procedures to compute the set of the Arf good semigroups of $\mathbb{N}^r$ with a fixed conductor $\textbf{c} \in \mathbb{N}^r$, for all $r \in \mathbb{N}$. 
\end{abstract}

\let\thefootnote\relax

\section*{Introduction}

The aim of this paper is to present some procedures regarding  Arf subsemigroups of $\mathbb{N}^r$ for $r \geq 1$. The concept of Arf numerical semigroup  firstly arises from the problem of studying the equivalence between two algebroid branches. Given an algebroid branch $R$, its multiplicity sequence is defined to be the sequence of the multiplicities of the succesive blowups $R_i$ of $R$. Two algebroid branches are equivalent if they have the same multiplicity sequence (cf. \cite[Definition 1.5.11]{Campi}). This approach, more geometric in nature, was translated in algebraic terms by Arf. In \cite{Arf} he introduced the concept of Arf ring and showed that for each ring $R$ there is a smallest Arf overring $R'$, called the Arf closure of $R$, that has also the same multiplicity sequence of $R$.  In the same paper it is proved that two algebroid branches are equivalent if and only if their Arf closure have the same value semigroup, that is a numerical Arf semigroup i.e. a numerical semigroup $S$ such that $S(s)-s$ is a semigroup, for each $s\in S$, where $S(s)=\left\{n \in S; n\geq s \right\}$.
In the numerical case an Arf semigroup $S=\left\{ s_0=0<s_1<s_2,\ldots \right\}$ is completely described by its multiplicity sequence, that is the sequence of the differences $s_{i+1}-s_i$.  So it is possible to study the Arf numerical semigroups by focusing on the properties of these sequences with a combinatorial approach. In this paper we give new procedures to calculate all Arf numerical semigroups with a prescribed conductor (this problem was already addressed in \cite{algo}). Then we address the Arf semigroups of  $\mathbb{N}^r$ focusing on the class of good semigroups. The concept of good semigroup was introduced in \cite{BDF}. Its definition depends on the properties of the value semigroup of a one-dimensional analytically unramified ring (for example the local ring of an algebraic curve), but in the same paper it is shown that the class of good semigroups is  bigger than the class of value semigroups. Therefore, the good semigroups can be seen as a natural generalization of the numerical semigroups and can be studied without referring to the ring theory context, with a more combinatorical approach. In this paper we deal only with local good semigroups.  A good semigroup $S \subseteq \mathbb{N}^n$ is local if the zero vector is the only element of $S$ with zero coordinates.  
We extend the definition  of Arf semigroup to the good semigroups of $\mathbb{N}^r$ in a natural way considering the usual partial ordering of $\mathbb{N}^r$, where $ \textbf{v}\leq \textbf{u}$ if and only if $v[i]\leq w[i]$ for all $i=1,\ldots,r$.  Extending the concept of multiplicity sequence, in \cite{BDF} it is also shown that to each local Arf good semigroup can be associated a multiplicity tree that characterizes the semigroup completely. A tree $T$ of vectors of $\mathbb{N}^n$ has to satisfy some properties to be a multiplicity tree of a local Arf good semigroup. For instance it must have multiplicity sequences along its branches (because the projection are Arf numerical semigroups) and each node must be able to be expressed as a sum of nodes in a subtree of $T$ rooted in it.
In \cite{Zito}, taking in account this $1$-$1$ correspondence, it is presented a  way to describe a multiplicity tree of an Arf semigroup of  $\mathbb{N}^r$  by mean of an ordered collection $E$ of multiplicity sequences and a matrix $M(T)_E=(p_{i,j})$. In the same paper it is shown that we can consider a particular class of trees, the untwisted ones, that are easier to study since they can be expressed by a vector of $\mathbb{N}^{r-1}$ instead of a matrix. This is not restrictive because every multiplicity tree can be transformed by an opportune permutation into an untwisted one (this implies that the corresponding semigroups are isomorphic).
 Then, using the properties of this representation and the procedure developed in the numerical case, we address the problem of computing the set of all the Arf good semigroups of $\mathbb{N}^r$ with an untwisted tree and a given conductor $\textbf{c}$, where the conductor $\textbf{c}$ of a good semigroup $S \subseteq \mathbb{N}^r$ is the minimal vector such that $ \textbf{c}+ \mathbb{N}^r \subseteq S$ (the existence of such a vector is guaranteed by the properties of good semigroups). 
 
The structure of the paper is the following.

In Section \ref{section1}, we firstly recall some definitions and properties concerning the numerical semigroups. Then we address the problem, already studied in \cite{algo}, of finding  the set of the multiplicity sequences of all the Arf numerical semigroups with a fixed conductor. In \cite{algo} the authors found a recursive algorithm for the computation of such a set, while in this section it is presented a non recursive procedure to determine it, that is faster than the previous one,  when used for large value of the conductor.

In Section \ref{section2}, we begin to deal with semigroups in higher dimensions. To do that we recall the definition of good semigroups and the properties of the multiplcity tree of an Arf local good semigroup. We denote by  $\textrm{Cond}(\textbf{c})$ the set of all the untwisted multiplicity trees of Arf good semigroups with conductor $\textbf{c} \in \mathbb{N}^r$. In order to find a recursive procedure for the computation of $\textrm{Cond}(\textbf{c})$, we start by solving the base case $r=2$.

In Section \ref{section3} we address the general case and, using Lemma \ref{Lemma3} and the base cases for $r=1$ and $r=2$, we are able to present a procedure that builds inductively the sets $\textrm{Cond}(\textbf{c})$ in all dimensions and for any value of the vector $\textbf{c}$.  We give a strategy for computing the set $\overline{\textrm{Cond}(\textbf{c})}$ of all the possible multiplicity trees (twisted and untwisted) associated to an Arf semigroup with conductor $\textbf{c}$. At the end of the section we give an example with $r=3$ on the computation of this set and we present some tables containing the cardinalities of the constructed sets for particular values of the conductor $\textbf{c}$.

The procedures presented here have been implemented in GAP (\cite{gap}).

\section{An algorithm for $\textrm{Cond}(n)$ where $n\in \mathbb{N}$} \label{section1}
In \cite{algo} it is presented an algorithm for the computation of the set of the Arf numerical semigroups with a given conductor.
In this section we give a new procedure for the computation of such a set that appeared to be faster when implemented in GAP.

First of all we recall some definitions and we fix some notations.
A numerical semigroup $S$ is a submonoid of $(\mathbb{N},+)$ such that $ \mathbb{N} \setminus S$ is finite. The maximum $F(S)$ of the set $ \mathbb{Z} \setminus S$ is known as the Frobenius number of $S$. The conductor $c(S)$ of $S$ is the smallest number such that $n \in S$ for all $n \geq c(S)$, and it is clear that we have $c(S)=F(S)+1$.

We are interested in a particular class of numerical semigroups, i.e. the set of Arf numerical semigroups. A numerical semigroup $S$ is said to be Arf if $S(s)-s$ is a semigroup for each $s \in S$, where $S(s)=\left\{ n \in S; n \geq s \right\}$. A nonincreasing sequence $M=\left\{m_n; n\geq 1 \right\}$ of $\mathbb{N}^*=\mathbb{N} \setminus \left\{ 0\right\}$  is called a multiplicity sequence if
\\
$ \bullet\ $ There exists $k \in \mathbb{N}^*$ such that $m_n=1$ for all $n \geq k$;
\\
$ \bullet\ $ For all $n \in \mathbb{N}^*$ there exists $s(n) \geq n+1$ such that $  m_n=\sum_{k=n+1}^{s(n)}{m_k}$.\\
In the following we  describe a multiplicity sequence $M$ by the finite vector $ M=[m_1,\ldots,m_k]$ where $m_k$ is the last entry in $M$ different from one. We  make an exception for the constant multiplicity sequence $M=\left\{1,\ldots,1,\ldots\right\}$ that will be represented by the vector $M=[1]$. 
Given a vector $\textbf{v}$ we always denote by $\textbf{v}[i]$ its $i$-th component

In \cite[Corollary 39]{Rosales} it is proved that a nonempty subset of $\mathbb{N}$ is an Arf numerical semigroup if and only if there exists a multiplicity sequence $M=[m_1,\ldots,m_k]$ such that \\ $S=\left\{0,m_1,m_1+m_2,\ldots,m_1+\dots+m_k,\rightarrow \right\},$
where with $n \rightarrow$ we  mean that all integers larger than $n$ are in $S$.
Thus we have a $1$-$1$ correspondence between Arf numerical semigroups and multiplicity sequences and, in the following, we will deal mainly with the latter.
If $M=[m_1,\ldots,m_k]$ is a multiplicity sequence, we call  $k$ the length of  $M$ and we will denote it  by $l(M)$.
We also denote by $\textrm{AS}(M)$ the Arf numerical semigroup associated to $M$. Now, from the previous result, it is clear that if $M$ is a mutliplicity sequence different from $[1]$, then the conductor of  $\textrm{AS}(M)$  is
$  \sum_{i=1}^{l(M)}{m_i},$
while the conductor of $\textrm{AS}( [1])=\mathbb{N}$ is $0$.
We denote by $\textrm{Cond}(n)$ the set of  the multiplicity sequences of  Arf numerical semigroups with conductor $n$.

We want to compute $\textrm{Cond}(n)$  with $ n \in \mathbb{N} $.
 If $n=0$ then  $\textrm{Cond}(n)=\left\{[1] \right\}$, while if   $n=1$ then  $\textrm{Cond}(n)=\emptyset.$ Thus, we suppose $n>1$.
Denote by $T^n(i)=\left\{ M \in \textrm{Cond}(i): M[1]+i \leq n \right\}  $ \textrm{ for all }$ i= \nolinebreak 2,\ldots,n-2.$
Now suppose that $M=[m_1,\ldots,m_k] \in \textrm{Cond}(n)$. If $k=1$ then $M=[n]$, otherwise we have the following situation:
\\
$ \bullet\ $ $
 2 \leq m_1<n-1 \textrm{ and } [m_2,\ldots,m_k] \in \textrm{Cond}(n-m_1)$;\\
 $ \bullet\ $ $
  m_1 \in \textrm{AS}([m_2,\ldots,m_k])$;\\
  $ \bullet\ $ $
  m_2-m_1 \leq 0 \Rightarrow m_2+n-m_1 \leq n \Rightarrow [m_2,\ldots,m_k] \in  T^n(n-m_1).$
  
Hence, if we know  $T^n(i)$  for $i=2,\ldots,n-2$, then we can compute $\textrm{Cond}(n)$ in the following way:
$$ \textrm{Cond}(n)=\bigcup_{i=2}^{n-2}{\left\{[n-i, M]: M \in T^n(i), n-i \in \textrm{AS}(M)\right\}} \cup \left\{  [n] \right\}.  $$
Now we need a way to compute  $T^n(i)$.
Suppose that $M=[m_1,\ldots,m_k] \in T^n(i)$. If $k=1$, and $2 \cdot i \leq n$ then $M=[i]$, otherwise we have the following situation:\\
$ \bullet\ 
 2 \leq m_1<i-1 \textrm{ and } [m_2,\ldots,m_k] \in \textrm{Cond}(i-m_1)=\textrm{Cond}(q);$  \\
$ \bullet\  q+m_2=i-m_1+m_2 \leq i \leq n \Rightarrow [m_2,\ldots,m_k] \in T^n(q);$ \\
  $ \bullet\ m_1 \in \textrm{AS}([m_2,\ldots,m_k]);$ \\
$\bullet\ m_1+i \leq n \Rightarrow 2m_1 \leq n+m_1-i \Rightarrow 2m_1 \leq n-q \Rightarrow m_1 \leq \left \lfloor{\frac{n-q}{2}}\right \rfloor.
$

So each $ T^n(i)$ can be constructed using $T^n(q)$ with $q<i$.
Thus, we have the following algorithm for  $\textrm{Cond}(n)$ for $n>1$.
\begin{algorithm}
\SetKwData{Left}{left}
\SetKwData{This}{this}
\SetKwData{Up}{up}
\SetKwFunction{Union}{Union}
\SetKwFunction{FindCompress}{FindCompress}
\SetKwInOut{Input}{input}
\SetKwInOut{Output}{output}
\caption{}
\Input{An integer $n>1$}
\Output{The set $\textrm{Cond}(n)$ of all the multiplicity sequences of Arf semigroups with conductor $n$}
\BlankLine
$\textrm{Cond}(n) \longleftarrow \left\{[n] \right\}$

\For{$i\leftarrow 2$ \KwTo $n-2$} {\If{$i\leq  \left \lfloor{\frac{n}{2}}\right \rfloor$}{$T^n(i) \longleftarrow \left\{ [i] \right\}$}\Else{$T^n(i) \longleftarrow \emptyset$}}
\For{$i\leftarrow 2$ \KwTo $n-2$}{
\For{$M \in T^n(i)$}{
\If{$n-i \in \textrm{AS}(M)$}{
$\textrm{Cond}(n) \longleftarrow \textrm{Cond}(n) \cup \left\{ [n-i,M]\right\}$
}
\For{$k \in \textrm{AS}(M) \cap \left\{ 2,\ldots,\left \lfloor{\frac{n-i}{2}}\right \rfloor \right\}$} {$T^n(i+k) \longleftarrow T^n(i+k) \cup \left\{ [k,M]\right\}$}
}
}
$\textrm{Cond}(n)$
\label{algo_disjdecomp}
\end{algorithm}
\newpage
\section{Arf good semigroups of $\mathbb{N}^2$ with given conductor} \label{section2}
From this section we begin to deal with Arf good semigroups of $\mathbb{N}^r$.  First of all we recall some important definitions and results.
 A good semigroup $S$ of $\mathbb{N}^r$ is a submonoid of $\left( \mathbb{N}^r,+\right)$ such that  (cf. \cite{BDF})\begin{itemize}
\item For all $a,b \in S$, $	\min(a,b) \in S$; 
\item If $a,b\in S$ and $a[i]=b[i]$ for some $i \in \left\{ 1, \ldots,n\right\}$, then there exists $c \in S$ such that $c[i] > a[i]=b[i]$, $c[j]\geq \min(a[j],b[j])$ for $j \in  \left\{ 1, \ldots,n\right\} \setminus \left\{ i\right\}$ and $c[j]=\min(a[j],b[j])$ if $a[j]\neq b[j]$;
\item There exists $\delta \in S$ such that $\delta+\mathbb{N}^r \subseteq S$
\end{itemize}
(where we are considering the usual partial ordering in $\mathbb{N}^r$: $a \leq b$ if $a[i] \leq b[i]$ for each $i= \nolinebreak 1,\ldots,r$). The previous properties guarantee that for a good semigroup $S$ there exists a minimal vector $c(S)$ such that $c(S)+ \mathbb{N}^r \subseteq S$. Such a vector is said to be the conductor of $S$.
\\
In this paper we  always deal with local good semigroups. A good semigroup $S$  is local if the zero vector is the only vector of $S$ with some component equal to zero.
However, it can be shown that every good semigroup is the direct product of local semigroups (cf. \cite[Theorem 2.5]{BDF}).

An Arf semigroup of $\mathbb{N}^r$, is a good semigroup such that $ S(\alpha)-\alpha $ is a semigroup, for each $\alpha \in S$ , where $ S(\alpha)=\left\{\beta \in S; \beta \geq \alpha \right\}$.
The  multiplicity tree $T$ of a local Arf semigroup $S\subseteq \mathbb{N}^r$ is a tree whose nodes are vectors $ \textbf{n}_i^j \in \mathbb{N}^r$, where with $\textbf{n}_i^j$ we mean that  this node is in the $i$-th branch on the $j$-th level (the root of the tree is $\textbf{n}_{1}^1=\textbf{n}_i^1 $ for all $i$ because we are in the local case and at level one all the branches must be glued) and we have 

$$ S=\left\{\textbf{0}\right\} \bigcup_{T'} \left\{ \sum_{\textbf{n}_i^j \in T' } {\textbf{n}_i^j}\right\},$$

where $T'$ ranges over all finite subtree of $T$ rooted in   $\textbf{n}_1^1$.

Furthermore a tree $T$ is a multiplicity tree of an Arf semigroup if and only if its  nodes  satisfy the following properties (cf. \cite[Theorem 5.11]{BDF} ).

\begin{itemize}
\item There exists $L \in \mathbb{N}$ such that for $m \geq L$, $\textbf{n}_i^m=(0,\ldots,0,1,0\ldots,0)$ (the nonzero coordinate is in the $i$-th position) for any $i=1,\ldots,n$;
\item $\textbf{n}_i^j[h]=0$ if and only if $\textbf{n}_i^j$ is not in the $h$-th branch of the tree;
\item Each $\textbf{n}_i^j$ can be obtained as a sum of nodes in a finite subtree $T'$ of $T$ rooted in $\textbf{n}_i^j$.
\end{itemize}
Notice that from these properties it follows that we must have multiplicity sequences along each branch. Therefore a multiplicity tree $T$ of an Arf semigroup of $\mathbb{N}^r$ can be represented by an ordered collection of $r$ multiplicity sequences $E$ and by an upper triangular matrix $ r \times r$  $$M(T)_{E}=\left(  \begin{matrix} 0 & p_{1,2} & p_{1,3} & \ldots & p_{1,r} \\ 
0 & 0 & p_{2,3} & \ldots & p_{2,r} \\ \ldots & \ldots & \ldots & \ldots & \ldots \\ 0 & 0 & 0 & \ldots& p_{r-1,r} \\ 0 & 0 & 0 & \ldots & 0 \\  \end{matrix}\right)$$where $p_{i,j}$ is the highest level such that the  $i$-th and the $j$-th branches are glued in $T$. Suppose that $E=\left\{M_1,\ldots,M_r \right\}$, then in \cite[Proposition 1.2]{Zito}, it is shown that the values that can be assigned to the $p_{i,j}$, in order to have a tree compatible with the previous properties, depend only on $M_i$ and $M_{j}$. In particular, if $M_i=M_{j},$ $p_{i,j}$ can be any positive integer,  while, if $M_i \neq M_{j}$, there exists a level $k_E(i,j) \in \mathbb{N}^*$ such that $p_{i,j}$ is at most $k_E(i,j)$.
If we set $s_{i,k}$ as the integer such that  $$ M_i[k]=\sum_{l=k+1}^{s_{i,k}}{M_i[l]},$$
then we have
 $$k_E(i,j)=\min\left\{ \min( s_{i,k}, s_{j,k}): s_{i,k}\neq s_{j,k} \right\}. $$ 
Thus, if $M_1$ and $M_2$  are two distinct multiplicity sequences, we  define the compatibility between $M_1$ and $M_2$ as the integer  $\textrm{Comp}(M_1,M_2)=k_E(1,2) $ where $E=\left\{M_1,M_2 \right\}$. While if $M_1= \nolinebreak M_2$ we set by definition  $\textrm{Comp}(M_1,M_2)=+\infty$. 
A tree $T$ is untwisted if two non-consecutive branches are glued at level $l$ if and only if  all the consecutive branches  between them are glued  at a level greater or equal to $l$. We call twisted a tree that  is not untwisted.
From the definition it follows that the matrix of an untwisted tree $T \in \tau(E)$ is such that: $$ p_{i,j}=\min\left\{ p_{i,i+1},\ldots,p_{j-1,j}\right\} \textrm{ for all } i<j. $$
So an untwisted tree can be completely described by the second diagonal of its matrix. Thus in the following we will  indicate an untwisted tree by a vector $T_E=(p_1,\ldots,p_{r-1})$ where $p_i=p_{i,i+1}$.
\begin{remark0} \label{remark0}
 It is easy to see that a twisted tree can be converted to an untwisted one by accordingly permuting its branches (the corresponding Arf semigroups are therefore isomorphic). 
Thus in the following we can focus only on the properties of the untwisted trees, that are easier to study, because the twisted ones can be always obtained by  permutation of an untwisted tree. 
\end{remark0}
The aim of this and the following section is to find a procedure that let us to determine all the local Arf semigroup  $S \subseteq\mathbb{N}^r$  with a given conductor $\textbf{c} \in \mathbb{N}^r$. 
For the Remark \ref{remark0}, we can focus only on the untwisted trees.

We denote by $\textrm{Cond}(\textbf{c})$  the set of all the untwisted multiplicity trees of  Arf semigroups in $\mathbb{N}^r$ with conductor $\textbf{c}\in \mathbb{N}^r$ (in the case $r=1$ we have the multiplicity sequences and from the previous section we have a procedure to determine such a set). 

We notice the following general fact.
\begin{prop} \label{prop}
Let $S$ be an Arf semigroup of $\mathbb{N}^r$, $T$ the corresponding multiplicity tree and $M_i$  for $i=1,\ldots,n $ the multiplicity sequences of its branches.

We introduce the following integers
$$ d(i)=\min\left\{  j \in \mathbb{N}: M_i[j]=1 \textrm{ and  the $i-$th branch is not glued to other branches at level $j$} \right\},$$ for $  i=1,\ldots,r .$

Then $\textbf{c}=(c[1],\ldots,c[n])$ is the conductor of $S$ where $$ c[i]=\sum_{k=1}^{d(i)-1}{M_i[k]}  \textrm{ for }  i=1,\ldots,n. $$
\end{prop}
\noindent \textbf{Proof}. Denote by $N(T)= \left\{ \textbf{n}_i^j \right\}$ the set of the nodes of $T$. We call $e_i=(0,\ldots,0,1,0,\ldots,0)$, where the non zero coordinate is in the $i$-th position. Now, from the definition of the integers $d(i)$, it follows that
\begin{itemize}
\item  $\textbf{n}_i^{d(i)}=e_i$ for all $i=1,\ldots,r$;
\item  $\textbf{n}_i^{d(i)-1}\neq e_i$ for all $i=1,\ldots,r$.
\end{itemize}
We consider the subtree $T'$ of $T$ such that $N(T')=\left\{\textbf{n}_i^{j(i)}: i=1,\ldots,r; j(i)=1,\ldots,d(i)-1 \right\}$.  Then  we have 
\begin{itemize}
\item $T'$ is rooted in  $\textbf{n}_1^{1}$ (it corresponds to an element of the associated Arf good semigroup);
\item $e_i \notin N(T')$ for all $i=1,\ldots,r$;
\item If $T''$ is such that $T' \subseteq T'' \subseteq T$ then $N(T'') \setminus N(T') $ consists only of nodes of the type $e_i$.
\end{itemize}
From the previous properties it is clear that the element corresponding to the subtree $T'$ must be the conductor of the Arf semigroup associated to $T$.
It is also trivial that the sum of all the elements of $N(T')$ is equal to 
$\textbf{c}=(c[1],\ldots,c[n])$ where $$ c[i]=\sum_{k=1}^{d(i)-1}{M_i[k]}  \textrm{ for }  i=1,\ldots,n. $$ \qed

\begin{remark1}
Using the notations of Proposition \ref{prop}, given an untwisted tree $T_{E}=(p_1,\ldots,p_{r-1})$, where $E=\left\{M_1,\ldots,M_r\right\}$, it is easy to show that \begin{flalign*}d(i)-1&=\max(l(M_i),p_i,p_{i-1}) \textrm{ for  } i=2,\ldots,r-1 \\ d(1)-1&=\max(l(M_1),p_1) \textrm { and } d(r)-1=\max(l(M_r),p_{r-1}). \end{flalign*}
\end{remark1}
Now, we focus on the case $r=2$ and we determine a procedure to compute $\textrm{Cond}(\textbf{c})$  where $\textbf{c}$ is a fixed arbitrary vector $(c[1],c[2])$. 
Suppose that $T_{E}=(p) \in \textrm{Cond}(\textbf{c})$ where $E=\left\{ M_1,M_2 \right\}$.
From the previous remark  $ d(1)-1=\max(l(M_1),p)$ and $d(2)-1=\max(l(M_2),p)$.
We have the following cases:
\\
$\bullet\ \textrm{Case}\  d(1)-1=l(M_1)$ and $d(2)-1=l(M_2)$.

We have $p \leq \min(l(M_1),l(M_2))$. Furthermore, we have $p \leq \textrm{Comp}(M_1,M_2)$  because $T$ is well defined.
Because $T_{E}=(p) \in \textrm{Cond}(\textbf{c})$ we have:
$$ c[1]=\sum_{k=1}^{d(1)-1}{M_1[k]} =\sum_{k=1}^{l(M_1)}{M_1[k]}  \textrm{ and }  c[2]=\sum_{k=1}^{d(2)-1}{M_2[k]} =\sum_{k=1}^{l(M_2)}{M_2[k]},$$
and  from it we deduce, if $c[1] \neq 1$ and $c[2] \neq 1 $, that $M_1 \in \textrm{Cond}(c[1])$  and  $M_2 \in \textrm{Cond}(c[2])$. If $c[1]=1$ then we can deduce $M_1=[1] \in \textrm{Cond}(0)$ (the same for $c[2]$). 
\begin{remark2}
 The multiplicity sequence $M=[1]$ is the only one such that the sum of its entries up to to its length is not equal to the conductor of its associated Arf numerical semigroup (that is $\mathbb{N}$ with conductor $0$).
In order to have a more compact notation in the following, we define the set $\textrm{Cond}^{*}(c)$ for all $c \in \mathbb{N} \setminus \left\{ 0 \right\}$, where
\begin{itemize}
 \item $\textrm{Cond}^{*}(c)=\textrm{Cond}(c)  \textrm{ for all } c \neq 1;$
\item $\textrm{Cond}^{*}(1)=\textrm{Cond}(0).$ 
\end{itemize}
\end{remark2}
So in this case $T$ belongs to the following set:
\begin{eqnarray*}  S^1(\textbf{c})=\left\{ T_{E}=(k): E=\left\{M_1,M_2\right\}; M_i \in\textrm{Cond}^{*}(c[i]) \textrm{  and } \right. \\  \left. 1 \leq  k \leq \min(\textrm{Comp}(M_1,M_2),l(M_1),l(M_2)) \right\}.  \end{eqnarray*}
On the other hand, we can notice that $ S^1(\textbf{c}) \subseteq  \textrm{Cond}(\textbf{c})$ (using the inverse implications).
\\
$\bullet\ \textrm{Case}\  d(1)-1=l(M_1)$ and $d(2)-1 \neq l(M_2)$.

Hence  $d(2)-1=p$  and $l(M_2)<p \leq \min(l(M_1),\textrm{Comp}(M_1,M_2)).$
Therefore:
$$ c[1]= \sum_{k=1}^{l(M_1)}{M_1[k]} \textrm{ , }  c[2]= \sum_{k=1}^{p}{M_2[k]}=\sum_{k=1}^{l(M_2)}{M_2[k]}+ \sum_{k=l(M_2)+1}^{p}{M_2[k]}=\sum_{k=1}^{l(M_2)}{M_2[k]}+p-l(M_2),$$
and from this we can deduce  $M_1 \in \textrm{Cond}^*(c[1])$  and, denoted by $k(2)=c[2]-(p-l(M_2))$, $M_2 \in \textrm{Cond}^*(k(2))$.
Notice that $k(2)<c[2]$.

Now, for all $k<c[2]$ we define the set:
\begin{eqnarray*} I_1(k)=\left\{ T_{E}=(p): E=\left\{M_1,M_2\right\}, M_1 \in \textrm{Cond}^*(c[1]),M_2 \in \textrm{Cond}^*(k) \textrm{ and }\right.\\ \left. p= l(M_2)+c[2]-k \leq \min(l(M_1),\textrm{Comp}(M_1,M_2)) \right\}. \end{eqnarray*}
Thus  $T$ belongs to
$$ S_2^1(\textrm{c})=\bigcup_{k=1}^{c[2]-1} {I_1(k)}.$$
With the inverse implication we can easily show that $ S_2^1(\textbf{c}) \subseteq  \textrm{Cond}(\textbf{c})$.
\\
$\bullet\ \textrm{Case}\  d(1)-1\neq l(M_1)$ e $d(2)-1 = l(M_2)$.

We have $d(1)-1=p$  and  $l(M_1)<p \leq \min(l(M_2),\textrm{Comp}(M_1,M_2)).$
Hence:
$$c[1]= \sum_{k=1}^{p}{M_1[k]}=\sum_{k=1}^{l(M_1)}{M_1[k]}+ \sum_{k=l(M_1)+1}^{p}{M_1[k]}=\sum_{k=1}^{l(M_1)}{M_1[k]}+p-l(M_1)  \textrm{; }  c[2]= \sum_{k=1}^{l(M_2)}{M_2[k]},$$
and from this we obtain $M_2 \in \textrm{Cond}^*(c[2])$ and, denoted by $k(1)=c[1]-(p-l(M_1))$, we deduce $M_1 \in \textrm{Cond}^*(k(1))$.
Notice that $k(1)<c[1]$.

For all $k<c[1]$ we define the set:
\begin{eqnarray*} I_2(k)=\left\{ T_{E}=(p): E=\left\{M_1,M_2\right\}, M_1 \in \textrm{Cond}^*(k),M_2 \in \textrm{Cond}^*(c[2]) \right. \textrm{ and }\\ \left. p= l(M_1)+c[1]-k \leq \min(l(M_2),\textrm{Comp}(M_1,M_2)) \right\}. \end{eqnarray*}
Therefore $T$ belongs to
$$ S_2^2(\textrm{c})=\bigcup_{k=1}^{c[1]-1} {I_2(k)}.$$
With the inverse implication we can easily show that  $ S_2^2(\textbf{c}) \subseteq  \textrm{Cond}(\textbf{c})$.
\\
$\bullet\ \textrm{Case}\  d(1)-1\neq l(M_1)$ and $d(2)-1 \neq l(M_2)$.
 
Then $d(2)-1=p$,  $d(1)-1=p$  and we have $\max(l(M_1),l(M_2))<p \leq \textrm{Comp}(M_1,M_2).$
It follows:
$$ c[1]= \sum_{k=1}^{p}{M_1[k]}=\sum_{k=1}^{l(M_1)}{M_1[k]}+ \sum_{k=l(M_1)+1}^{p}{M_1[k]}=\sum_{k=1}^{l(M_1)}{M_1[k]}+p-l(M_1) $$ $$  c[2]= \sum_{k=1}^{p}{M_2[k]}=\sum_{k=1}^{l(M_2)}{M_2[k]}+ \sum_{k=l(M_2)+1}^{p}{M_2[k]}=\sum_{k=1}^{l(M_2)}{M_2[k]}+p-l(M_2).$$
If we denote by $k(1)=c[1]-(p-l(M_1))$ and by $k(2)=c[2]-(p-l(M_2))$, we have $M_1 \in \nolinebreak \textrm{Cond}^*(k(1))$ and $M_2 \in \textrm{Cond}^*(k(2))$.

Furthermore, notice that $k(1)<c[1]$ and $k(2)<c[2]$.
Now, for all $k_1<c[1]$  and  $k_2<c[2]$ we define the set:
\begin{eqnarray*} I(k_1,k_2)=\left\{ T_{E}=(p): E=\left\{M_1,M_2\right\}, M_i \in \textrm{Cond}^*(k_i) \textrm{ for } i=1,2 \textrm{ and }\right.\\ \left.  p=l(M_1)+c[1]-k_1=l(M_2)+c[2]-k_2 \leq \textrm{Comp}(M_1,M_2) \right\}. \end{eqnarray*} So  $T$ belongs to $$   S^3(\textbf{c})=\bigcup_{1\leq k_i < c[i]}{I(k_1,k_2)} \subseteq \textrm{Cond}(\textbf{c}).$$
Even in this case we can show that $ S^3(\textbf{c}) \subseteq  \textrm{Cond}(\textbf{c})$.

We have studied all the possible cases so we proved 
$$ S^1(\textbf{c}) \cup S_1^2(\textbf{c}) \cup S_2^2(\textbf{c}) \cup S^3(\textbf{c})=\textrm{Cond}(\textbf{c}).$$
All the previous set can be computed by using the procedure given in the case $r=1$ so we have found a procedure to compute $\textrm{Cond}(\textbf{c})$ when $r=2$.
\begin{example1}
Let us compute $ \textrm{Cond}([4,5])$. 

First of all we compute $S^1([4,5])$.  We need $\textrm{Cond}^*(4)$ and $\textrm{Cond}^*(5)$. They are:
$$  \textrm{Cond}^*(4)=\left\{[4],[2,2]\right\} \textrm{ and } \textrm{Cond}^*(5)=\left\{ [5], [3,2] \right\}.$$
Hence when we compute $S^1([4,5])$ we find:
\begin{itemize}
\item $E_1=\left\{M_1=[4] , M_2=[5]\right\}$. Thus $\textrm{Comp}(M_1,M_2)=5$ and \newline
 $\min(l(M_1),l(M_2))=1$. Then we have only the tree $T_1=T_{E_1}=(1)$.
\item $E_2=\left\{M_1=[4] , M_2=[3,2]\right\}$. Thus $\textrm{Comp}(M_1,M_2)=3$ and  \newline
 $\min(l(M_1),l(M_2))=1$. Then we have only the tree $T_2=T_{E_2}=(1)$.
\item  $E_3=\left\{M_1=[2,2] , M_2=[5]\right\}$. Thus $\textrm{Comp}(M_1,M_2)=2$  and \newline
 $\min(l(M_1),l(M_2))=1$. Then we have only the tree $T_3=T_{E_3}=(1)$.
\item  $E_4=\left\{M_1=[2,2] , M_2=[3,2]\right\}$. We have $\textrm{Comp}(M_1,M_2)=2$ and \newline
 $\min(l(M_1),l(M_2))=2$. So we have the trees  $T_4=T_{E_4}=(1)$ and  $T_5=T_{E_4}=(2)$ .
\end{itemize}
Hence $S^1([4,5])=\left\{ T_1,T_2,T_3,T_4,T_5\right\}$.
\vskip 0.3in

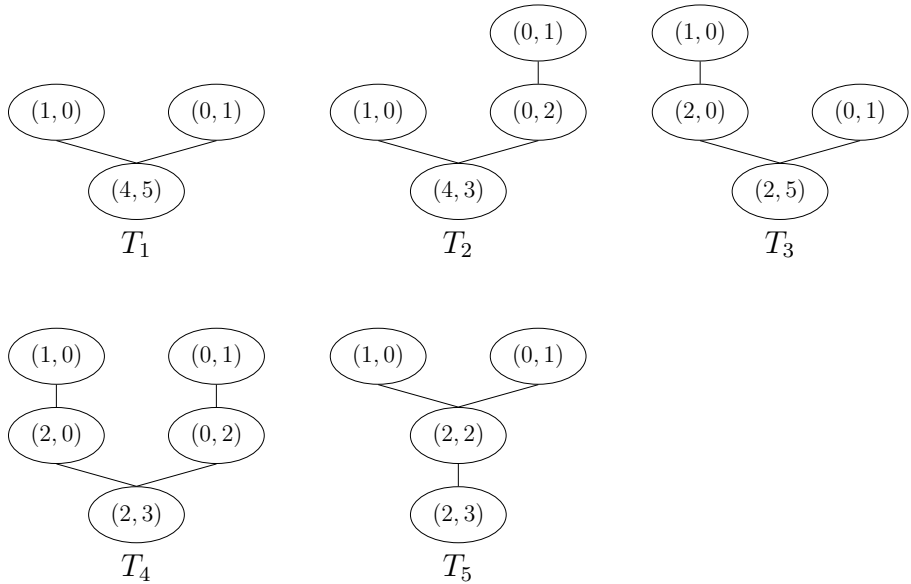
\begin{tikzpicture}[grow'=up,sibling distance=32pt,scale=.75]
\tikzset{level distance=40pt,every tree node/.style={draw,ellipse}} \Tree [ .$(4,5)$  [ .$(1,0)$  ] [ .$(0,1)$ ]  ] ]\node[below]at(current bounding box.south){$T_1$}; \end{tikzpicture}
\hskip 0.3in
\begin{tikzpicture}[grow'=up,sibling distance=32pt,scale=.75]
\tikzset{level distance=40pt,every tree node/.style={draw,ellipse}} \Tree [ .$(4,3)$  [ .$(1,0)$  ] [ .$(0,2)$  [.$(0,1)$ ]]  ]\node[below]at(current bounding box.south){$T_2$}; \end{tikzpicture}
\hskip 0.3in
\begin{tikzpicture}[grow'=up,sibling distance=32pt,scale=.75]
\tikzset{level distance=40pt,every tree node/.style={draw,ellipse}} \Tree [ .$(2,5)$  [ .$(2,0)$ [.$(1,0)$ ] ] [ .$(0,1)$ ]  ]\node[below]at(current bounding box.south){$T_3$}; \end{tikzpicture}
\vskip 0.3in
\begin{tikzpicture}[grow'=up,sibling distance=32pt,scale=.75]
\tikzset{level distance=40pt,every tree node/.style={draw,ellipse}} \Tree [ .$(2,3)$  [ .$(2,0)$ [.$(1,0)$ ] ] [ .$(0,2)$ [.$(0,1)$ ] ]  ]\node[below]at(current bounding box.south){$T_4$}; \end{tikzpicture}
\hskip 0.3in
\begin{tikzpicture}[grow'=up,sibling distance=32pt,scale=.75]
\tikzset{level distance=40pt,every tree node/.style={draw,ellipse}} \Tree [.$(2,3)$ [ .$(2,2)$  [ .$(1,0)$  ] [ .$(0,1)$ ]  ] ] ]\node[below]at(current bounding box.south){$T_5$}; \end{tikzpicture}

  Now we compute   $S_1^2([4,5])$. The only value $k$  such that  $I_1(k) \neq \emptyset $  is $k=4$ and we have:
\begin{itemize}
\item $c[1]=4$, $k=4$. If we consider $M_1=[2,2] \in \textrm{Cond}^*(4)$ , \newline $M_2=[4] \in \textrm{Cond}^*(4) $ and  $E_5=\left\{ M_1,M_2 \right\}$ we have $$ l(M_2)+c[2]-k=2 \leq \min(l(M_1),\textrm{Comp}(M_1,M_2))=\min(2,2)=2. $$ Hence we have the tree $T_6=T_{E_5}=(2).$ 
\end{itemize}
Therefore  $S_1^2([4,5])=\left\{T_6\right\}$.
Let us compute $S_2^2([4,5])$. The only value $k$  such that  $I_2(k) \neq \emptyset$  is $k=3$ and we have:
\begin{itemize}
\item $k=3$, $c[2]=5$.  If we consider $M_1=[3] \in \textrm{Cond}^*(3)$ , \newline $M_2=[3,2] \in \textrm{Cond}^*(5) $ and  $E_6=\left\{ M_1,M_2 \right\}$ we have $$ l(M_1)+c[1]-k=1+4-3=2 \leq \min(l(M_2), \textrm{Comp}(M_1,M_2)=\min(2,3)=2. $$
Hence we have the tree $T_7=T_{E_6}=(2).$ 
\end{itemize}
Therefore  $S_2^2([4,5])=\left\{T_7\right\}$.

We finally compute  $S^3([4,5])$.

 The only values of $k_1$  and $k_2$  such that  $I(k_1,k_2) \neq \emptyset$  are the following:
\begin{itemize}
\item $k_1=2$, $k_2=3$.  If we consider $M_1=[2] \in \textrm{Cond}^*(2)$, \newline $M_2=[3] \in \textrm{Cond}^*(3) $  and  $E_7=\left\{ M_1,M_2 \right\}$ we have $$ l(M_1)+c[1]-k_1=1+4-2=3=1+5-3=l(M_2)+c[2]-k_2 \leq \textrm{Comp}(M_1,M_2)=3. $$
 Thus we have the tree  $T_8=T_{E_7}=(3).$ 
\item $k_1=3$, $k_2=4$.  If we consider $M_1=[3] \in \textrm{Cond}^*(3)$ , \newline $M_2=[4] \in \textrm{Cond}^*(4) $ and $E_8=\left\{ M_1,M_2 \right\}$ we have $$ l(M_1)+c[1]-k_1=1+4-3=2=1+5-4=l(M_2)+c[2]-k_2\leq \textrm{Comp}(M_1,M_2)=4. $$
Thus we have the tree $T_9=T_{E_8}=(2).$ 
Hence $S^3([4,5])=\left\{ T_8,T_9\right\}$.
\vskip 0.3in

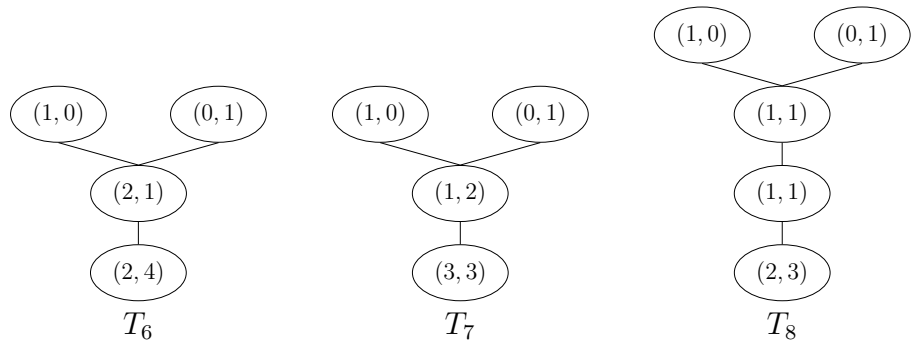
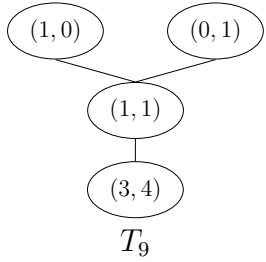
\begin{tikzpicture}[grow'=up,sibling distance=32pt,scale=.75]
\tikzset{level distance=40pt,every tree node/.style={draw,ellipse}} \Tree [.$(2,4)$ [ .$(2,1)$  [ .$(1,0)$  ] [ .$(0,1)$ ]  ] ] ]\node[below]at(current bounding box.south){$T_6$}; \end{tikzpicture}
\hskip 0.3in
\begin{tikzpicture}[grow'=up,sibling distance=32pt,scale=.75]
\tikzset{level distance=40pt,every tree node/.style={draw,ellipse}} \Tree [.$(3,3)$ [ .$(1,2)$  [ .$(1,0)$  ] [ .$(0,1)$ ]  ] ] ]\node[below]at(current bounding box.south){$T_7$}; \end{tikzpicture}
\hskip 0.3in
\begin{tikzpicture}[grow'=up,sibling distance=32pt,scale=.75]
\tikzset{level distance=40pt,every tree node/.style={draw,ellipse}} \Tree [.$(2,3)$ [.$(1,1)$ [ .$(1,1)$  [ .$(1,0)$  ] [ .$(0,1)$ ]  ] ] ]]\node[below]at(current bounding box.south){$T_8$}; \end{tikzpicture}
\vskip 0.3in
\begin{tikzpicture}[grow'=up,sibling distance=32pt,scale=.75]
\tikzset{level distance=40pt,every tree node/.style={draw,ellipse}} \Tree [.$(3,4)$ [ .$(1,1)$  [ .$(1,0)$  ] [ .$(0,1)$ ]  ] ] ]\node[below]at(current bounding box.south){$T_9$}; \end{tikzpicture}
\end{itemize}
Summarizing, we have $ \textrm{Cond}([4,5])=\left\{ T_1,T_2,T_3,T_4,T_5, T_6,T_7,T_8,T_9\right\}$.
\end{example1}
\begin{example2}
Using the previous results it is easy to implement an algorithm that computes the number of Arf semigroups of $ \mathbb{N}^2$ with a given conductor. Each entry of the following table is such a number, where the conductors range from $(1,1)$ to  $(20,20)$.

\begin{center}
	\resizebox{16.5cm}{4.5cm}{
	\begin{tabular}{| | l || c| c | c | c| c| c| c| c| c| c| c| c| c| c| c|c| c| c| c| c ||p{5cm} }
		\hline \hline 
		 & \textbf{1} &\textbf{2}  & \textbf{3} &\textbf{4} & \textbf{5}&\textbf{6}&\textbf{7}&\textbf{8}&\textbf{9}&\textbf{10}&\textbf{11}&\textbf{12}&\textbf{13}&\textbf{14}&\textbf{15}&\textbf{16}&\textbf{17}&\textbf{18}&\textbf{19}&\textbf{20} \\ \hline \hline
		\textbf{1} & 1 & 1  &1 &2 & 2&4&3&7&6&10&9&17&12&25&20&32&27&49&34&68  \\ \hline
		\textbf{2} & 1 & 2  &2 &4 &4&8&6&14&12&20&18&34&24&50&40&64&54&98&68&136  \\  \hline
		\textbf{3} &1& 2& 3& 4& 5& 9& 7& 16& 14& 22& 21& 39& 26& 57& 46& 71& 60& 111& 75& 155  \\
		\hline
		\textbf{4} &2& 4& 4& 10& 9& 18& 15& 33& 28& 49& 43& 81& 59& 120& 96& 156& 131& 236& 167& 328 \\
		\hline
		\textbf{5} & 2& 4& 5& 9& 12& 19& 15& 34& 32& 51& 45& 86& 62& 128& 102& 161& 139& 250& 172& 347  \\
		\hline
		\textbf{6}  & 4& 8& 9& 18& 19& 41& 30& 68& 60& 99& 92& 171& 122& 252& 201& 326& 275& 497& 344& 687   \\
		\hline
		\textbf{7}  & 3& 6& 7& 15& 15& 30& 30& 54& 48& 80& 74& 134& 104& 204& 163& 264& 221& 399& 285& 556   \\
		\hline
	\textbf{8}  &  7& 14& 16& 33& 34& 68& 54& 129& 108& 180& 164& 306& 222& 453& 371& 593& 499& 901& 632& 1251\\
	\hline
	
		\textbf{9}  &  6& 12& 14& 28& 32& 60& 48& 108& 108& 160& 147& 271& 202& 404& 330& 522& 459& 809& 566& 1120 \\
		\hline
		\textbf{10}  &   10& 20& 22& 49& 51& 99& 80& 180& 160& 284& 242& 454& 337& 676& 545& 878& 748& 1336& 961& 1867\\
		\hline
		\textbf{11}  &    9& 18& 21& 43& 45& 92& 74& 164& 147& 242& 245& 412& 307& 611& 502& 798& 685& 1215& 868& 1688 \\
		\hline
		\textbf{12}  &    17& 34& 39& 81& 86& 171& 134& 306& 271& 454& 412& 798& 567& 1148& 927& 1492& 1273& 2277& 1608& 3159 \\
		\hline
		\textbf{13}  &     12& 24& 26& 59& 62& 122& 104& 222& 202& 337& 307& 567& 469& 849& 694& 1115& 961& 1689& 1224& 2347 \\
		\hline
		\textbf{14}  &      25& 50& 57& 120& 128& 252& 204& 453& 404& 676& 611& 1148& 849& 1750& 1383& 2224& 1897& 3389& 2403& 4710\\
		\hline
		\textbf{15}  &      20& 40& 46& 96& 102& 201& 163& 371& 330& 545& 502& 927& 694& 1383& 1192& 1805& 1556& 2753& 1976& 3822\\
		\hline
		\textbf{16}  &       32& 64& 71& 156& 161& 326& 264& 593& 522& 878& 798& 1492& 1115& 2224& 1805& 2992& 2493& 4433& 3174& 6155 \\
		\hline
		\textbf{17}  &        27& 54& 60& 131& 139& 275& 221& 499& 459& 748& 685& 1273& 961& 1897& 1556& 2493& 2244& 3798& 2734& 5266 \\
		\hline
		\textbf{18}  &        49& 98& 111& 236& 250& 497& 399& 901& 809& 1336& 1215& 2277& 1689& 3389& 2753& 4433& 3798& 6867& 4814& 9394 \\
		\hline
		\textbf{19}  &        34& 68& 75& 167& 172& 344& 285& 632& 566& 961& 868& 1608& 1224& 2403& 1976& 3174& 2734& 4814& 3634& 6701 \\
		\hline
		\textbf{20}  &        68& 136& 155& 328& 347& 687& 556& 1251& 1120& 1867& 1688& 3159& 2347& 4710& 3822& 6155& 5266& 9394& 6701& 13219 \\
		\hline \hline

	\end{tabular}}
\end{center}
\end{example2}
\section{Arf semigroups of  $\mathbb{N}^r$ with a given conductor} \label{section3}
In this section we study the general case. We want to develope a recursive procedure to calculate $\textrm{Cond}(\textbf{c})$ for $\textbf{c} \in \mathbb{N}^r$, using the fact that we already know how to solve the base cases $r=1$ and $r=2$. In order to do that is very useful the following Lemma.

\begin{Lemma3} \label{Lemma3}
Consider $\textbf{c}=(c[1],\ldots,c[r]) \in \mathbb{N}^r$, with $r \geq 3 $  and suppose that the untwisted tree $T=T_{E}=(p_1,\ldots,p_{r-1}) \in \textrm{Cond}(\textbf{c})$, where $E=\left\{ M_1,\ldots,M_r \right\}$. If $t \in \left\{ 2,\ldots,r-1\right\}$, then we have that at least one of the following conditions must hold: \begin{itemize}
\item $ T_1= T_{E_1}=(p_1,\ldots,p_{t-1}) \in \textrm{Cond}((c[1],\ldots,c[t]))$  with $E_1=\left\{M_1,\ldots,M_t\right\};$
\item $T_2= T_{E_2}=(p_t,\ldots,p_{r-1}) \in \textrm{Cond}((c[t],\ldots,c[r]))$  with $E_2=\left\{M_t,\ldots,M_r\right\}.$
\end{itemize}
\end{Lemma3}
\noindent \textbf{Proof}.
We assume by contradiction that
 \begin{itemize}
\item $ T_1 \notin\textrm{Cond}((c[1],\ldots,c[t])); $
\item $ T_2 \notin \textrm{Cond}((c[t],\ldots,c[r])). $
\end{itemize}
Let us consider the following integers (which are clearly linked to the conductor): \\
$ d(i)=\min\left\{  j \in \mathbb{N}: M_i[j]=1 \textrm{ and the $i$-th branch in $T$ is not glued to other branches at level $j$} \right\}$, for $  i=1,\ldots,r .$\\
$ d_1(i)=\min\left\{  j \in \mathbb{N}: M_i[j]=1 \textrm{  and the $i$-th branch in $T_1$ isn't glued to other branches at level $j$}  \right\}$, for $  i=1,\ldots,t .$ \\
$ d_2(i)=\min\left\{  j \in \mathbb{N}: M_i[j]=1 \textrm{ and the $i$-th branch in $T_2$ isn't glued to other branches at level $j$}  \right\}$, for $  i=t,\ldots,r .$\\
We  have $d_1(l)=d(l)$, for all $l=1,\ldots,t-1$,  and  $d_2(m)=d(m)$,  for all $m=t+1,\ldots,r$. Furthermore, $d_1(t)\leq d(t)$ and $d_2(t) \leq d(t)$. In fact
we have noticed that $d(t)-1=\max(l(M_t),p_{t-1},p_t)$, while $d_1(t)-1=\max(l(M_t),p_{t-1})$  and  $d_2(t)-1=\max(l(M_t),p_t)$.

From  $T \in \textrm{Cond}(\textbf{c})$  we deduce that
 $$ \sum_{k=1}^{d(i)-1}{M_i[k]}=c[i]  \textrm{ for }  i=1,\ldots,r.$$
We denote respectively by $(c_1[1],\ldots,c_1[t])$  and by $(c_2[t],\ldots,c_2[r]))$ the conductors of $T_1$ and  $T_2$. 

We have:
 $$ c_1[l]=\sum_{k=1}^{d_1(l)-1}{M_l[k]}=\sum_{k=1}^{d(l)-1}{M_l[k]}=c[l] \textrm{   for  } l=1,\ldots,t-1  $$ $$ \textrm{  and } c_2[m]=\sum_{k=1}^{d_2(m)-1}{M_m[k]}=\sum_{k=1}^{d(m)-1}{M_m[k]}=c[m],  \textrm{ for } m=t+1,\ldots,r$$
and this implies, because $ T_1 \notin \textrm{Cond}((c[1],\ldots,c[t])) $ and  $ T_2 \notin \textrm{Cond}((c[t],\ldots,c[r]))$, that $$  c_1[t]=\sum_{k=1}^{d_1(t)-1}{M_t[k]}\neq c[t] \textrm{   and  }  c_2[t]=\sum_{k=1}^{d_2(t)-1}{M_t[k]} \neq c[t],$$
and therefore we have $d_1(t)<d(t)$ and $d_2(t)<d(t)$.

From this it would follow \begin{flalign*}
d(t)-1&=\max(l(M_t),p_{t-1},p_t)=\max(\max(l(M_t),p_{t-1}),\max(l(M_t),p_t)))=\\&=\max(d_1(t)-1,d_2(t)-1)<d(t)-1
\end{flalign*} and we obtain a contradiction. \qed

Now, using this Lemma, we can introduce an algorithm that solves our problem working inductively.
Given $\textbf{c} \in \mathbb{N}^r$, with $r \geq 3$,  we want to compute $\textrm{Cond}(\textbf{c})$. We suppose that we are able to solve the problem for all $s<r$ and we develope a strategy for the $r$ case.

Let us fix some notations.  If $k=2,\ldots,r-1$, we denote by $\textbf{c}_k=(c[1],\ldots,c[k])$ and by  $\textbf{c}^k=(c[k+1],\ldots,c[r])$.
Similarly, if $E=\left\{M_1,\ldots,M_r\right\}$, we denote by $E_k=\left\{M_1,\ldots,M_k\right\}$ and by  $E^k=\left\{M_{k+1},\ldots,M_r\right\}.$  Furthermore, for $i=1,\ldots,r-1$, we define the integers ${}^*p_i= \nolinebreak \max(l(M_i),p_{i-1})$ and $p_i^*=\max(l(M_{i+1}),p_{i+1})$, where, by definition, we set $p_{r-1}^*=l(M_r)$ and  ${}^*p_1=l(M_1)$.

Fixed $\textbf{c} \in \mathbb{N}^r$, we suppose to have a tree $T=T_{E}=(p_1,\ldots,p_{r-1}) \in \textrm{Cond}(\textbf{c})$ with $E=\nolinebreak \left\{ M_1,\ldots,M_r \right\}$.
Consider $t \in \left\{ 2,\ldots,r-1 \right\}$. It follows from Lemma \ref{Lemma3} that we  only have two cases:
\\
$\bullet\ \textrm{Case}\  T_1=T_{E_t}=(p_1,\ldots,p_{t-1}) \in \textrm{Cond}(\textbf{c}_t).$

We clearly have $d_1(i)=d(i) $ for all $i=1,\ldots,t-1$, while from $T_1 \in \textrm{Cond}(\textbf{c}_t)$ it follows that $d_1(t)=d(t)$.
Hence:
$$ ^*p_t=\max(l(M_t),p_{t-1})=d_1(t)-1=d(t)-1=\max(p_t,p_{t-1},l(M_t))=\max(^*p_t,p_t)$$
and we deduce that $p_t \leq  {}^*p_t$.

We consider the tree $T_2=T_{E^t}=(p_{t+1},\ldots,p_{r-1})$, (if $t=r-1$ we have $T_2=M_r$).
We clearly have $d_2(i)=d(i)$ for all $i=t+2,\ldots,r$.

On the other hand $d_2(t+1)-1=\max(l(M_{t+1}),p_{t+1})=p_{t}^*$  may be different from $$d(t+1)-1= \max(l(M_{t+1},p_{t+1},p_t)=\max(p_t^*,p_t).$$
Hence we have the following two subcases:
\\
$ \blacktriangleright \textrm{Subcase}\ d_2(t+1)=d(t+1).$

In this case we have $T_2 \in \textrm{Cond}(\textbf{c}^t)$ and  $p_t \leq p^*_t$ (if $t=r-1$ we have $T_2=M_r \in \textrm{Cond}^*(c[r])$). We also recall that we must have the compatibility condition $p_t \leq \textrm{Comp}(M_t,M_{t+1})$.

Thus we have discovered that $T$ belongs to the following set:
\begin{equation*}  S_1^1(\textbf{c})=\left\{ T_{E}=(p_1,\ldots,p_{r-1}): E=\left\{M_1,\ldots,M_r \right\}; T_{E_t}=(p_1,\ldots,p_{t-1})\in \textrm{Cond}(\textbf{c}_t);\right. \end{equation*}  \begin{equation*} \left. T_{E^t}=(p_{t+1},\ldots,p_{r-1}) \in \textrm{Cond}(\textbf{c}^t) \textrm{ with }1\leq p_t \leq \min({}^*p_{t},p^*_{t},\textrm{Comp}(M_t,M_{t+1})) \right\}.\end{equation*} 
It is very easy to check that we also have  $S_1^1(\textbf{c}) \subseteq \textrm{Cond}(\textbf{c}).$
If  $t=r-1$  the previous set has the following definition:
\begin{equation*}  S_1^1(\textbf{c})=\left\{ T_{E}=(p_1,\ldots,p_{r-1}): E=\left\{M_1,\ldots,M_r \right\}; \right.  T_{E_{r-1}}=(p_1,\ldots,p_{r-2}) \in \textrm{Cond}(\textbf{c}_{r-1});\end{equation*}  \begin{equation*} \left.  M_r \in \textrm{Cond}^*(c[r]) ;1\leq p_{r-1} \leq \min(^*p_{r-1},l(M_r),\textrm{Comp}(M_{r-1},M_r)) \right\}.\end{equation*}
\\
$ \blacktriangleright \textrm{Subcase}\ d_2(t+1) \neq d(t+1).$

In this case we have $$p_t=d(t+1)-1 >d_2(t+1)-1=\max(l(M_{t+1}),p_{t+1})=p^*_t. $$
Hence $$ c[t+1]=\sum_{k=1}^{d(t+1)-1} M_{t+1}[k]=\sum_{k=1}^{p_t} M_{t+1}[k] =\sum_{k=1}^{p^*_t} M_{t+1}[k]+p_t-p^*_t $$
and from this it  follows that $T_2 \in \textrm{Cond}((k[t+1],c[t+2],\ldots,c[r]))$, where $$\sum_{k=1}^{p^*_t} M_{t+1}[k]=\sum_{k=1}^{d_2(t+1)-1} M_{t+1}[k]=k[t+1]<c[t+1],$$
and we have $T_2 \in \textrm{Cond}^*(k[t+1])$ in the case $t=r-1$.
Thus we have $k[t+1]=c[t+1]-(p_t-p^*_t)$.
For all the $k_{t+1} \in \mathbb{N}$ such that $k_{t+1}<c[t+1]$ we define the set:
\begin{equation*} I_1(k_{t+1})=\left\{ T_{E}=(p_1,\ldots,p_{r-1}):E=\left\{M_1,\ldots,M_r \right\}; T_{E_t}=(p_1,\ldots,p_{t-1}) \in\textrm{Cond}(\textbf{c}_t); \right. \end{equation*}  \begin{equation*} \left.  T_{E^t}=(p_{t+1},\ldots,p_{r-1}) \in \textrm{Cond}((k_{t+1},c[t+2],\ldots,c[r])); \right.\end{equation*} \begin{equation*} \left.  p_t=p_t^* +c[t+1]-k_{t+1}\leq\min(^*p_{t},\textrm{Comp}(M_t,M_{t+1}))\right\}.\end{equation*}
Hence $T$ belongs to the following set:
$$ S_1^2(\textbf{c})= \bigcup_{k_{t+1}=1}^{c[t+1]-1}{I_1(k_{t+1})}, $$
and it is clear that $ S_1^2(\textbf{c}) \subseteq   \textrm{Cond}(\textbf{c})$.

If $t=r-1$ the previous set has the following definition:
\begin{equation*} I_1(k_r)=\left\{ T_{E}=(p_1,\ldots,p_{r-1}):E=\left\{M_1,\ldots,M_r \right\}; T_{E_{r-1}}=(p_1,\ldots,p_{r-2}) \in \textrm{Cond}(\textbf{c}_{r-1});\right. \end{equation*}  \begin{equation*} \left. M_r \in \textrm{Cond}^*(k_r) ;p_{r-1}=l(M_r)+c[r]-k_r\leq \min(^*p_{r-1},\textrm{Comp}(M_{r-1},M_r)) \right\}.\end{equation*}
\\
$\bullet\ \textrm{Case}\  T_2=T_{E^{t-1}}=(p_t,\ldots,p_{r-1})\in \textrm{Cond}(\textbf{c}^{t-1}).$

We  only have to adapt the considerations made in the previous case to this case.  Thus we directly give the sets which arise without further justifications.
 \\
 $ \blacktriangleright $ If  $t\neq 2$,
\begin{equation*}  S_2^1(\textbf{c})=\left\{ T_{E}=(p_1,\ldots,p_{r-1}): E=\left\{M_1,\ldots,M_r \right\}; T_{E^{t-1}}=(p_t,\ldots,p_{r-1})\in \textrm{Cond}(\textbf{c}^{t-1});\right. \end{equation*}  \begin{equation*} \left. T_{E_{t-1}}=(p_1,\ldots,p_{t-2}) \in \textrm{Cond}(\textbf{c}_{t-1}) \textrm{ con }1\leq p_{t-1} \leq \min(^*p_{t-1},p^*_{t-1},\textrm{Comp}(M_t,M_{t-1})) \right\}.\end{equation*} 
$ \blacktriangleright $  If $t=2$,
\begin{equation*}  S_2^1(\textbf{c})=\left\{ T_{E}=(p_1,\ldots,p_{r-1}): E=\left\{M_1,\ldots,M_r \right\}; \right.  T_{E^1}=(p_2,\ldots,p_{r-1}) \in \textrm{Cond}(\textbf{c}^{1});\end{equation*}  \begin{equation*} \left.  M_1 \in \textrm{Cond}^*(c[1]) ;1\leq p_{1} \leq \min(l(M_1),p^*_{1},\textrm{Comp}(M_{1},M_2)) \right\}.\end{equation*}
We have $S_2^1(\textbf{c}) \subseteq  \textrm{Cond}(\textbf{c}).$

For all $k_{t-1} \in \mathbb{N}$ such that $1\leq k_{t-1} < c[t-1]$ we consider:
 \\
$ \blacktriangleright $ If  $t \neq 2,$ 
\begin{equation*} I_2(k_{t-1})=\left\{ T_{E}=(p_1,\ldots,p_{r-1}):E=\left\{M_1,\ldots,M_r \right\}; T_{E^{t-1}}=(p_t,\ldots,p_{r-1}) \in\textrm{Cond}(\textbf{c}^{t-1}); \right. \end{equation*}  \begin{equation*} \left.  T_{E_{t-1}}=(p_{1},\ldots,p_{t-2}) \in \textrm{Cond}((c[1],\ldots,c[t-2],k_{t-1}));\right.\end{equation*} \begin{equation*} \left. p_{t-1}={}^*p_{t-1}+c[t-1]-k_{t-1}\leq  \min(p^*_{t-1},\textrm{Comp}(M_{t-1},M_{t}))\right\}.\end{equation*}
$ \blacktriangleright $  If $t=2$,
\begin{equation*} I_2(k_1)=\left\{ T_{E}=(p_1,\ldots,p_{r-1}):E=\left\{M_1,\ldots,M_r \right\}; T_{E^{1}}=(p_2,\ldots,p_{r-1}) \in \textrm{Cond}(\textbf{c}^{1})\right. \end{equation*}  \begin{equation*} \left. M_1 \in \textrm{Cond}^*(k_1) ;p_{1}=l(M_1)+c[1]-k_1\leq \min(p^*_{1},\textrm{Comp}(M_{1},M_2)) \right\}.\end{equation*}
We have that:
$$ S_2^2(\textbf{c})= \bigcup_{ k_{t-1}=1}^{c[t-1]-1}{I_2(k_{t-1})} \subseteq   \textrm{Cond}(\textbf{c}).$$
The previous lemma ensures that we have considered all the possibilities. 
So we showed that 
$$\textrm{Cond}(\textbf{c}) \subseteq  S_1^1(\textbf{c}) \cup S_1^2(\textbf{c}) \cup S_2^1(\textbf{c}) \cup S_2^2(\textbf{c}), $$
hence 
$$ S_1^1(\textbf{c}) \cup S_1^2(\textbf{c}) \cup S_2^1(\textbf{c}) \cup S_2^2(\textbf{c})=\textrm{Cond}(\textbf{c}).$$
Due to our induction hypothesis  all the previous sets can be computed  so we developed an algorithm which computes $\textrm{Cond}(\textbf{c})$.

Now we have a way to compute all the untwisted multiplicity trees with a given conductor $\textbf{c}$ for all the $\textbf{c} \in \mathbb{N}^r$. Suppose that we want to find also the twisted multiplicity trees with conductor $\textbf{c}$.
We will call $\overline{\textrm{Cond}(\textbf{c})}$ the set of all multiplicity trees (twisted or untwisted) associated to an Arf semigroup with conductor $\textbf{c}$.
Suppose that $T$ is a twisted tree in $\overline{\textrm{Cond}(\textbf{c})}$ with $\textbf{c} \in \mathbb{N}^r$. Then there exists a permutation $\sigma \in S^r$, where $S^r$ is the symmetric group, such that $\sigma(T)$ is untwisted and it clearly belongs to $ \textrm{Cond}(\sigma(\textbf{c}))$. From this it follows that:
$$ \overline{\textrm{Cond}(\textbf{c})}=\bigcup_{\sigma \in S^r} \left\{ \sigma^{-1}(T): T \in  \textrm{Cond}(\sigma(\textbf{c})) \right\}.$$
\begin{example3}
Let us compute $ \textrm{Cond}([3,2,4])$. In this case $r=3$, therefore we have $t=2$.
First  of all we compute $S_1^1([3,2,4])$. Because $t=r-1$ the definition of this set is:
\begin{equation*}  S_1^1([3,2,4])=\left\{ T_{E}=(p_1,p_{2}): E=\left\{M_1,M_2,M_3 \right\}; \right.  T_{E_{2}}=(p_1) \in \textrm{Cond}([3,2]);\end{equation*}  \begin{equation*} \left.  M_3 \in \textrm{Cond}^*(4) ;1\leq p_{2} \leq \min(l(M_3),^*p_{2},\textrm{Comp}(M_{2},M_3)) \right\}.\end{equation*}
Then to do that we need the follwing sets: 
\begin{itemize} 
\item $\textrm{Cond}([3,2])= \left\{A_1,A_2 \right\}$ where $A_1=T_{F_1}=(1)$ and $A_2=T_{F_2}=(2)$ with $F_1=\left\{ [3],[2]\right\}$ and $F_2=\left\{ [2],[1]\right\}.$
\item  $\textrm{Cond}^*(4)=\left\{[2,2],[4] \right\}$
\end{itemize}
Hence we consider: \begin{itemize} \item $E_1=\left\{ M_1=[3],M_2=[2],M_3=[2,2] \right\}$  and we have \newline $\min(\max(l(M_2),p_1),\textrm{Comp}(M_2,M_3),l(M_3))=\min(1,2,2)=1$. Thus we only have  the tree $T_1=T_{E_1}=(1,1).$
\item $E_2=\left\{ M_1=[3],M_2=[2],M_3=[4] \right\}$ and we have \newline $\min(\max(l(M_2),p_1),\textrm{Comp}(M_2,M_3),l(M_3))=\min(1,3,1)=1$. Thus we only have  the tree $T_2=T_{E_2}=(1,1).$
\item $E_3=\left\{ M_1=[2],M_2=[1],M_3=[2,2] \right\}$  and we have \newline $\min(\max(l(M_2),p_1),\textrm{Comp}(M_2,M_3),l(M_3))=\min(2,2,2)=2$. Thus we have the trees $T_3=T_{E_3}=(2,1)$ and  $T_4=T_{E_3}=(2,2).$
\item  $E_4=\left\{ M_1=[2],M_2=[1],M_3=[4] \right\}$ and we have \newline $\min(\max(l(M_2),p_1),\textrm{Comp}(M_2,M_3),l(M_3))=\min(2,2,1)=1$. Thus we only have  the tree $T_5=T_{E_4}=(2,1).$
\end{itemize}
Hence $S_1^1([3,2,4])=\left\{ T_1,T_2,T_3,T_4,T_5\right\}.$

Now we compute $S_1^2([3,2,4])$. We find $k_3=3$ as the only value such that $I(k_3) \neq \emptyset$. In fact, if we consider  $A_2$ and $M_3=[3]$, we have:
\begin{itemize}
\item $ E_5=\left\{ M_1=[2],M_2=[1],M_3=[3] \right\} \textrm{ and we have } $ 

$
 l(M_3)+c[3]-k_3=2 \leq \min(\max(l(M_2),p_1),\textrm{Comp}(M_2,M_3)))=\min(2,2)=2$. Thus we have the tree $T_6=T_{E_5}=(2,2).$
  \end{itemize}
\vskip 0.3in

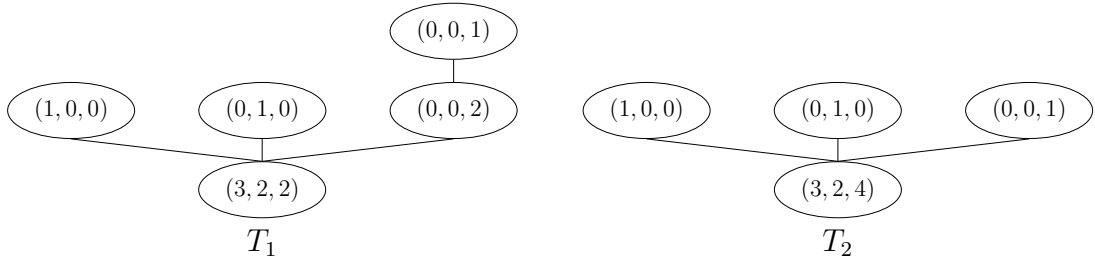
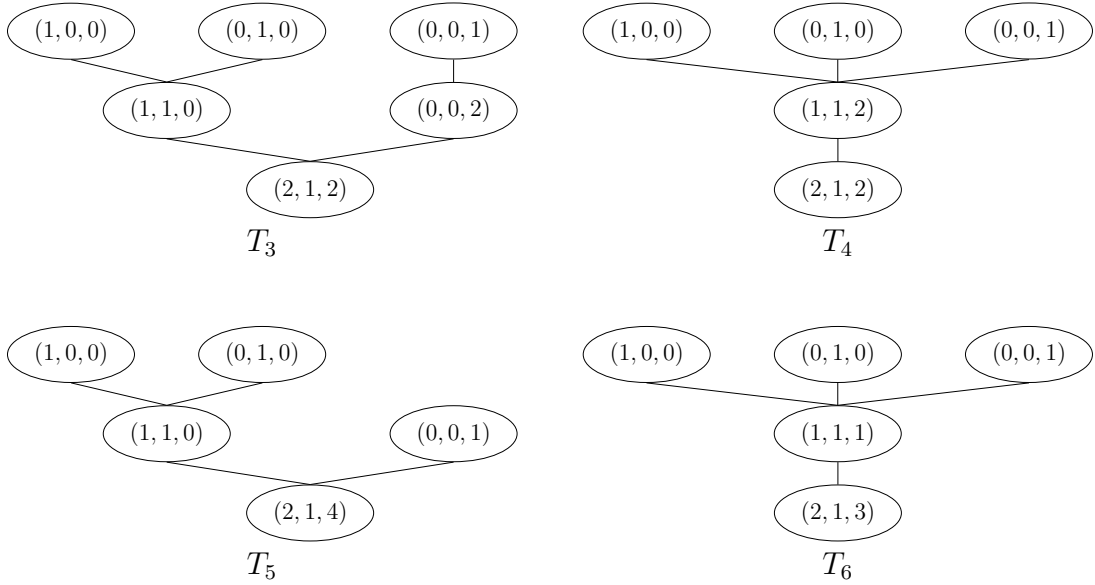
\begin{tikzpicture}[grow'=up,sibling distance=32pt,scale=.75]
\tikzset{level distance=40pt,every tree node/.style={draw,ellipse}} \Tree [ .$(3,2,2)$ [ .$(1,0,0)$ ]  [ .$(0,1,0)$ ] [ .$(0,0,2)$ [ .$(0,0,1)$ ] ] ] \node[below]at(current bounding box.south){$T_1$}; \end{tikzpicture}
\hskip 0.3in
\begin{tikzpicture}[grow'=up,sibling distance=32pt,scale=.75]
\tikzset{level distance=40pt,every tree node/.style={draw,ellipse}} \Tree [ .$(3,2,4)$ [ .$(1,0,0)$ ]  [ .$(0,1,0)$ ] [ .$(0,0,1)$ ] ] \node[below]at(current bounding box.south){$T_2$}; \end{tikzpicture}
\vskip 0.3in
\begin{tikzpicture}[grow'=up,sibling distance=32pt,scale=.75]
\tikzset{level distance=40pt,every tree node/.style={draw,ellipse}} \Tree [ .$(2,1,2)$ [ .$(1,1,0)$   [  .$(1,0,0)$ ]  [ .$(0,1,0)$  ] ] [ .$(0,0,2)$ [ .$(0,0,1)$ ] ] ] \node[below]at(current bounding box.south){$T_3$}; \end{tikzpicture}
\hskip 0.3in
\begin{tikzpicture}[grow'=up,sibling distance=32pt,scale=.75]
\tikzset{level distance=40pt,every tree node/.style={draw,ellipse}} \Tree [ .$(2,1,2)$ [ .$(1,1,2)$ [ .$(1,0,0)$ ] [ .$(0,1,0)$ ] [ .$(0,0,1)$ ]  ] ] \node[below]at(current bounding box.south){$T_4$}; \end{tikzpicture}
\vskip 0.3in
\begin{tikzpicture}[grow'=up,sibling distance=32pt,scale=.75]
\tikzset{level distance=40pt,every tree node/.style={draw,ellipse}} \Tree [ .$(2,1,4)$ [ .$(1,1,0)$   [  .$(1,0,0)$ ]  [ .$(0,1,0)$  ] ] [ .$(0,0,1)$  ] ] \node[below]at(current bounding box.south){$T_5$}; \end{tikzpicture}
\hskip 0.3in
\begin{tikzpicture}[grow'=up,sibling distance=32pt,scale=.75]
\tikzset{level distance=40pt,every tree node/.style={draw,ellipse}} \Tree [ .$(2,1,3)$ [ .$(1,1,1)$ [ .$(1,0,0)$ ] [ .$(0,1,0)$ ] [ .$(0,0,1)$ ]  ] ] \node[below]at(current bounding box.south){$T_6$}; \end{tikzpicture}

Now we compute $S_2^1([3,2,4])$. We are in the case $t=2$ so its definition is:
\begin{equation*}  S_2^1([3,2,4])=\left\{ T_{E}=(p_1,p_2): E=\left\{M_1,M_2,M_3 \right\}; \right.  T_{E^1}=(p_2) \in \textrm{Cond}([2,4]);\end{equation*}  \begin{equation*} \left.  M_1 \in \textrm{Cond}^*(3) ;1\leq p_{1} \leq \min(l(M_1),p^*_{1},\textrm{Comp}(M_{1},M_2)) \right\}.\end{equation*}
Then to do that we need the following sets: 
$$\textrm{Cond}([2,4])= \left\{B_1,B_2,B_3,B_4 \right\}, \textrm{ where }$$ 
\begin{itemize} 
\item  $B_1=T_{G_1}=(2)$  with $G_1=\left\{ [1],[2,2]\right\};$ 
\item  $B_2=T_{G_2}=(1)$  with $G_2=\left\{ [2],[2,2]\right\}$; 
\item  $B_3=T_{G_3}=(2)$  with $G_3=\left\{ [1],[3]\right\}$ ;
\item  $B_4=T_{G_4}=(1)$  with $G_4=\left\{ [2],[4]\right\}$ ;
\item  $\textrm{Cond}^*(3)=\left\{[3] \right\}.$
\end{itemize}
Hence we consider: \begin{itemize} \item $E_6=\left\{ M_1=[3],M_2=[1],M_3=[2,2] \right\}$ and we have \newline $\min(\max(p_2,l(M_2)),\textrm{Comp}(M_1,M_2),l(M_1))=\min(2,2,1)=1$. Thus we only have  the tree $T_7=T_{E_6}=(1,2).$
\item $E_1=\left\{ M_1=[3],M_2=[2],M_3=[2,2] \right\}$ and we have \newline $\min(\max(p_2,l(M_2)),\textrm{Comp}(M_1,M_2),l(M_1))=\min(1,3,1)=1$.  Thus we only have  the tree,  already found  in   $S_1^1([3,2,4])$, $T_1=T_{E_1}=(1,1).$ 
\item $E_7=\left\{ M_1=[3],M_2=[1],M_3=[3] \right\}$ and we have \newline $\min(\max(p_2,l(M_2)),\textrm{Comp}(M_1,M_2),l(M_1))=\min(2,2,1)=1$. Hence we have the tree  $T_8=T_{E_7}=(1,2)$.
\item  $E_2=\left\{ M_1=[3],M_2=[2],M_3=[4] \right\}$ and we have \newline $\min(\max(p_2,l(M_2)),\textrm{Comp}(M_1,M_2),l(M_1))=\min(1,3,1)=1$. Thus we only have  the tree, already found in   $S_1^1([3,2,4])$, $T_2=T_{E_2}=(1,1).$
\end{itemize}
Hence  $S_2^1([3,2,4])=\left\{ T_1,T_2,T_7,T_8\right\}.$

Now we compute  $S_2^2([3,2,4])$. We find $k_1=2$ as the only value such that  $I(k_1) \neq \emptyset$, and $I(2)$ contains two elements. In fact, if we consider $B_1$ and $B_3$ and $M_1=[2]$, we have:
\begin{itemize}
\item $ E_3=\left\{ M_1=[2],M_2=[1],M_3=[2,2] \right\} \textrm{ and we have } $

$ l(M_1)+c[1]-k_1=2 \leq \min(\max(p_2,l(M_2)),\textrm{Comp}(M_1,M_2)))=\min(2,2)=2$. Thus we only have  the tree, already found in   $S_1^1([3,2,4])$,  $T_4=T_{E_3}=(2,2).$
\item $ E_5=\left\{ M_1=[2],M_2=[1],M_3=[3] \right\} \textrm{ and we have } $

$ l(M_1)+c[1]-k_1=2 \leq \min(\max(p_2,l(M_2)),\textrm{Comp}(M_1,M_2)))=\min(2,2)=2$. Thus we only have  the tree, already found in   $S_1^2([3,2,4])$,   $T_6=T_{E_5}=(2,2).$
  \end{itemize}
Thus $S_2^1([3,2,4]) \cup S_2^2([3,2,4])=\left\{ T_1,T_2,T_4,T_6,T_7,T_8\right\}$.
\vskip 0.3in
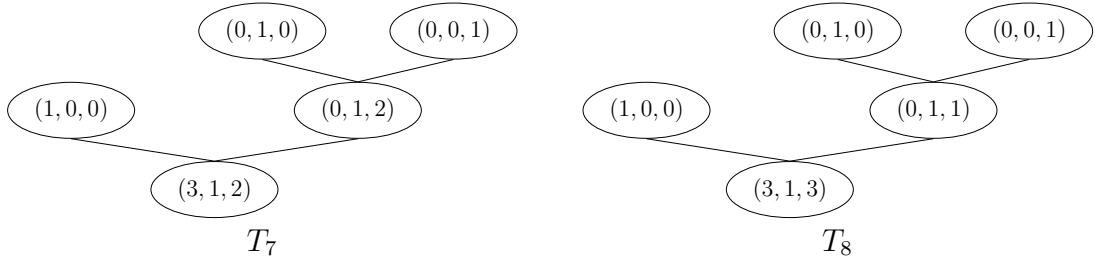
\begin{tikzpicture}[grow'=up,sibling distance=32pt,scale=.75]
\tikzset{level distance=40pt,every tree node/.style={draw,ellipse}} \Tree [ .$(3,1,2)$ [ .$(1,0,0)$ ] [ .$(0,1,2)$ [ .$(0,1,0)$ ] [ .$(0,0,1)$ ]  ]  ] \node[below]at(current bounding box.south){$T_7$}; \end{tikzpicture}
\hskip 0.3in
\begin{tikzpicture}[grow'=up,sibling distance=32pt,scale=.75]
\tikzset{level distance=40pt,every tree node/.style={draw,ellipse}} \Tree [ .$(3,1,3)$ [ .$(1,0,0)$ ] [ .$(0,1,1)$ [ .$(0,1,0)$ ] [ .$(0,0,1)$ ]  ]  ] \node[below]at(current bounding box.south){$T_8$}; \end{tikzpicture}

Hence,  $ \textrm{Cond}([3,2,4])=\left\{ T_1,T_2,T_3,T_4,T_5,T_6,T_7,T_8\right\}$.

If we compute the set $\overline{\textrm{Cond}([3,2,4])}$, with the technique explained above, we find that:
$$ \overline{\textrm{Cond}([3,2,4])}=\textrm{Cond}([3,2,4]) \bigcup \left\{ T_9, T_{10} \right\},$$
where 
\begin{itemize}
\item $T_9=M(T)_{E_8}=\left( \begin{matrix} 0 & 1 & 2 \\ 0 & 0 & 1 \\ 0 &0 & 0 \end{matrix}\right)$ where $E_8=\left\{ M_1=[2], M_2=[2], M_3=[2,2] \right\}$.
\item  $T_{10}=M(T)_{E_9}=\left( \begin{matrix} 0 & 1 & 2 \\ 0 & 0 & 1 \\ 0 &0 & 0 \end{matrix}\right)$ where $E_9=\left\{ M_1=[2], M_2=[2], M_3=[3] \right\}$.
\end{itemize}
\vskip 0.3in
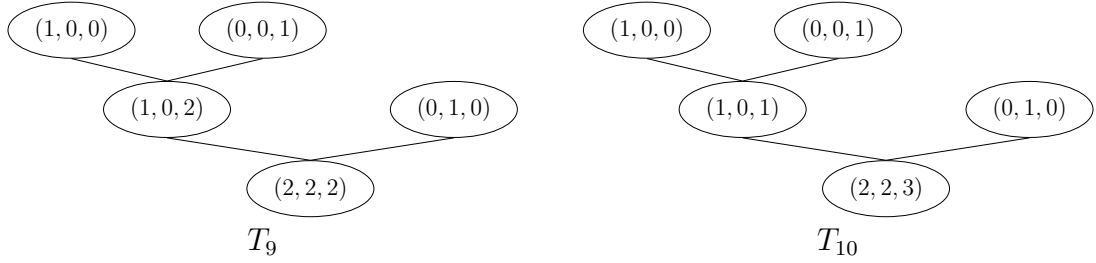
\begin{tikzpicture}[grow'=up,sibling distance=32pt,scale=.75]
\tikzset{level distance=40pt,every tree node/.style={draw,ellipse}} \Tree [ .$(2,2,2)$ [ .$(1,0,2)$   [  .$(1,0,0)$ ]  [ .$(0,0,1)$  ] ] [ .$(0,1,0)$  ] ] \node[below]at(current bounding box.south){$T_9$}; \end{tikzpicture}
\hskip 0.3in
\begin{tikzpicture}[grow'=up,sibling distance=32pt,scale=.75]
\tikzset{level distance=40pt,every tree node/.style={draw,ellipse}} \Tree [ .$(2,2,3)$ [ .$(1,0,1)$   [  .$(1,0,0)$ ]  [ .$(0,0,1)$  ] ] [ .$(0,1,0)$  ] ] \node[below]at(current bounding box.south){$T_{10}$}; \end{tikzpicture}

\end{example3}

\begin{example4}
It is easy to implement an algorithm that computes the number of untwisted  Arf semigroups of $ \mathbb{N}^3$ with a given conductor.  In the following table we have the values obtained for some conductors.
\begin{center}
		\resizebox{16cm}{!}{
	\begin{tabular}{ || c | c || c |c|| c| c|| c| c||}
		\hline \hline
	 $\textbf{c}$ & $| \textrm{Cond}(\textbf{c})|$& $\textbf{c}$ & $| \textrm{Cond}(\textbf{c})|$& $\textbf{c}$ & $| \textrm{Cond}(\textbf{c})|$& $\textbf{c}$ & $| \textrm{Cond}(\textbf{c})|$ \\ 
	 
	  \hline \hline
	 $[1,1,1]$ & 1 & $[8,8,8]$ & 2401 & $[15,15,15]$ & 71736& $[7,8,9]$ & 843 \\
	 \hline
	 $[2,2,2]$ & 4 & $[9,9,9]$ & 1940 & $[1,2,3]$ & 2& $[8,9,10]$ & 2901 \\
	 	 \hline
	 $[3,3,3]$ & 9 & $[10,10,10]$ & 8126 & $[2,3,4]$ & 8& $[9,10,11]$ & 3913 \\
	 	 \hline
	 $[4,4,4]$ & 50 & $[11,11,11]$ & 6671 & $[3,4,5]$ & 18& $[10,11,12]$ & 11178 \\
	 	 \hline
	 $[5,5,5]$ & 72 & $[12,12,12]$ & 37750 & $[4,5,6]$ & 86& $[11,12,13]$ & 13942 \\
	 	 \hline
	 $[6,6,6]$ & 425 & $[13,13,13]$ & 18263 & $[5,6,7]$ & 144& $[12,13,14]$ & 40278 \\
	 	 \hline
	 $[7,7,7]$ & 294 & $[14,14,14]$ & 123498 & $[6,7,8]$ & 542& $[13,14,15]$ & 47675 \\
	  \hline \hline
	\end{tabular}}
\end{center}
\end{example4} 
\begin{example44}
The following table contains the value of $|\overline{\textrm{Cond}(\textbf{c})}|$ for some values of $\textbf{c}$.
\begin{center}
	\resizebox{16cm}{!}{
		\begin{tabular}{ || c | c || c |c|| c| c|| c| c||}
			\hline \hline
			$\textbf{c}$ & $\big| \overline{\textrm{Cond}(\textbf{c})}\big| $& $\textbf{c}$ & $\big| \overline{\textrm{Cond}(\textbf{c})}\big| $& $\textbf{c}$ & $\big| \overline{\textrm{Cond}(\textbf{c})}\big| $& $\textbf{c}$ & $\big| \overline{\textrm{Cond}(\textbf{c})}\big| $ \\ 
			
			\hline \hline
			$[1,1,1]$ & 1 & $[7,7,7]$ & 406 & $[1,2,3]$ & 2& $[7,8,9]$ & 1145 \\
			\hline
			$[2,2,2]$ & 5 & $[8,8,8]$ &3217 & $[2,3,4]$ & 10& $[8,9,10]$ & 3828 \\
			\hline
			$[3,3,3]$ & 12 & $[9,9,9]$ & 2650& $[3,4,5]$ & 26& $[9,10,11]$ & 5289 \\
			\hline
			$[4,4,4]$ & 66 & $[10,10,10]$ & 10992 & $[4,5,6]$ & 110& $[10,11,12]$ & 14908 \\
			\hline
			$[5,5,5]$ & 98 & $[11,11,11]$ & 9131 & $[5,6,7]$ & 192& $[11,12,13]$ & 19147 \\
			\hline
			$[6,6,6]$ & 567 & $[12,12,12]$ & 50903 & $[6,7,8]$ & 701& $[12,13,14]$ & 53144\\
		
			\hline \hline
	\end{tabular}}
\end{center}
\end{example44}

 \begin{acknowledgements}
	The author would like to thank Marco D'Anna for his helpful comments and suggestions. Special thanks to  Pedro Garc\'ia-S\'anchez for his careful reading of an earlier version of the paper and for many helpful hints regarding the implementation in GAP of the presented procedures.
\end{acknowledgements}
\begin{otherlanguage}{english}

\end{otherlanguage}

GIUSEPPE ZITO-Dipartimento di Matematica e Informatica-Universit\`a di Catania-Viale Andrea Doria, 6, I-95125 Catania- Italy.

E-mail address: giuseppezito@hotmail.it


\begin{thebibliography}{Dillo 100}
\bibitem{Arf} C. Arf {\em Une interpretation algebrique de la suite des ordres de multiplicite d'une branche algebrique}  Proc. London Math. Soc. (2), 50:256-287, 1948.
\bibitem{BDF}  V.~Barucci, M.~D'Anna, R.~Fr\"oberg {\em Analitically unramifed one-dimensional semilocal rings and their value semigroups}, in J. Pure Appl. Algebra,  {\em 147} (2000),  215-254.
\bibitem{Campi} A. Campillo,  {\em Algebroid curves in positive characteristic},Lecture Notes in Math. 813, Springer-Verlag, Heidelberg, (1980).
\bibitem{algo} P. A. Garcia-Sanchez, B. A. Heredia, H. I. Karakas,  J. C. Rosales  {\em Parametrizing Arf numerical semigroups}, J. Algebra Appl. Vol. 16, No. 11 (2017) .
\bibitem{Rosales} J. C. Rosales, { \em Principal ideals of numerical semigroups}, Bull. Belg. Math. Soc. Simon Stevin 10 (2003), 329 - 343.
\bibitem{gap} The GAP Group,  {\em GAP-Groups, Algorithms and Programming}, Version 4.7.5 , {\em 25} (2014).
\bibitem{Zito} G. Zito, {\em Arf good semigroups}, to appear in Journal of Algebra and its Applications (2018).
\end{thebibliography}
\end{document}